\documentclass[12pt]{article}
\usepackage{amsfonts,amssymb,amscd,amsmath,amsthm}
\def\R{{\bf R}}   
\long\def\comment#1\endcomment{}

\def\sign{\mathop{\fam0 sign}}

\def\id{\mathop{\fam0 id}}

\def\im{\mathop{\fam0 im}}

\def\eps{\varepsilon}
\usepackage{graphicx}
\def\mpfile#1#2{\includegraphics{#1#2.eps}}

\begin{document}

\centerline{\uppercase{\bf Basic embeddings and Hilbert's 13th problem}
\footnote{
This is an English version of the paper in Russian under the same title. 
The English version has much shorter first section (which corresponds to two 
sections in Russian version), but contains solutions of problems 14a and 16c 
from the third section. 
Whenever possible I give references to surveys not to original papers.
I would like to acknowledge V.I.Arnold, Yu.M. Burman, I.N. Shnurnikov, A.R. 
Safin, S.M. Voronin and M. Vyaliy for useful discussions, and M. Vyaliy for 
preparation of figures.}
}
\smallskip
\centerline{\bf A. Skopenkov
\footnote{skopenko@mccme.ru, 
http://dfgm.math.msu.su/people/skopenkov/papersc.ps}
}

\bigskip
\small
{\bf Abstract.}
This note is purely expository. 
In the course of the Kolmogorov-Arnold solution of Hilbert's 13th problem on 
superpositions there appeared the notion of {\it basic embedding}.  
A subset $K$ of $\R^2$ is {\it basic} if for each continuous 
function $f\colon K\to\R$ there exist continuous functions $g,h\colon\R\to\R$ 
such that $f(x,y)=g(x)+h(y)$ for each point $(x,y)\in K$.
%We show how in the course of the Kolmogorov-Arnold solution of Hilbert's 13th 
%problem on superpositions there appeared the notion of a basic embedding. 
%A subset  K  of  R^2  is {\it basic} if for each continuous function  f:K->R  
%there exist continuous functions  g,h:R->R  such that  f(x,y) = g(x) + h(y)  
%for each point  (x,y)  in  K.
We present descriptions of basic subsets of the plane (with a proof) and 
description of graphs basically embeddable into the plane (solutions of 
Arnold's and Sternfeld's problems). 
We present some results and open problems on the smooth version of the property 
of being basic.  
This note is accessible to undergraduates and could be an interesting easy 
reading for mature mathematicians. 
The two sections can be read independently on each other. 

\normalsize

%Russian version has extended first two sections, but does not contain solutions 
%of problems 14a and 16c from the third section. 
%The main result of this sequence of problems (problem 8b) is an elementary
%solution [MT03] of the 'half' of the Arnold problem [Ar58'] on
%characterization of basic subsets of the plane [St89].
%See also [Vo81, Vo82, Sk95, Ku00, Ku03].
%The most important unsolved problems here are on characterization of
%{\it smoothly basic} subsets of the plane [RZ06].

\bigskip
\centerline{\uppercase{\bf Hilbert's 13th problem and basic embeddings}}

\smallskip
{\bf Hilbert's 13th problem} 

Let us recall informally the concept of {\it superposition}.
Suppose that there is a set of functions of several variables, 
including all variables considered as functions. 
Represent each of the functions as an element of a circuit with several entries 
and one exit.
Then a {\it superposition} of functions of this set is a function that can be 
repsesented by a circuit constructed from given elements; the circuit should 
not contain oriented cycles.

For example, a polynomial $a_n x^n+a_{n-1} x^{n-1}+\dots+a_1 x+a_0$ is a 
superposition of the constant functions and the functions $f(x,y)=x+y$, 
$g(x,y)=xy$.
It is clear that any elementary function can be represented as a
superposition  of functions of at most two variables.
{\it Is it possible to represent each function of several arguments as a 
superposition of functions of at most two arguments?}

Since there is a 1--1 correspondence between a segment and a square, any 
function of three and more variables is superposition of (in general, 
discontinuous) functions of two variables.
So the above question is only interesting for continuous functions. 
\footnote{
Denote by 
$$|x,y|=|(x_1,\dots,x_n),\ (y_1,\dots,y_n)|=
\sqrt{(x_1-y_1)^2+\dots+(x_n-y_n)^2}$$
the ordinary distance between points $x=(x_1,\dots,x_n)$ and 
$y=(y_1,\dots,y_n)$ of $\R^n$.
Let $K$ be a subset of $\R^n$.
A function $f:K\to\R$ is called {\it continuous} if for each point
$x_0\in K$ and number $\eps >0$  there exists a number $\delta>0$ such that
for each point $x\in K$ if $|x,x_0|<\delta$, then $|f(x)-f(x_0)|<\eps$.
E. g. the function $f(x_1,x_2)=\sqrt{x_1^2+x_2^2}$ is continuous on the plane, 
whereas the function $f(x_1,x_2)$ equal to the integer part of $x_1+x_2$ is 
not.}
From now on we assume all functions to be continuous, unless the contrary 
is explicitly specified. 

%We can set $\delta=\eps$. Then the statement follows from the triangle 
%inequality: $|f(z)-f(z_0)|\le |z,z_0|$.
%For the point $(1,0)$ and $\eps= \frac 12$ there is no such $\delta$ because
% $|f(1,0)-f(1-\frac\delta 2,0)|=1>\frac 12$.
%For example, any function on a finite set is continuous.
%Continuity of a function defined on a sequence of points $a_n\in R^2$
%converging to a point $a_0$ means that
%$\lim\limits_{n\to\infty} f(a_n)$ exists and equals $f(a_0)$.

\smallskip
{\bf Hilbert's 13th problem.}
{\it Can the equation $x^7+ax^3+bx^2+cx+1=0$ of degree seven be solved without 
using functions of three variables?}

\smallskip
This question was answered affirmatively in 1957 by Kolmogorov and Arnold.
They proved that any continuous function of $n$ variables defined on a compact 
subset of $\R^n$ can be represented as a superposition of continuous functions 
of one variable and addition. 
For an exposition accessible to undergraduates see~[Ar58]. 
See also [Vi04].

\bigskip
{\bf Basic embeddings into higher-dimensional spaces
\footnote{This subsection is not used in the sequel and so can be omitted.}
}

Ostrand extended the Kolmogorov-Arnold Theorem this theorem to arbitrary 
$n$-dimensional compacta [St89]. 
It is in the Kolmogorov-Arnold-Ostrand papers that the notion of basic 
subset appeared for the first time.
It was explicitly introduced by Sternfeld [St89]. 
A subset $K\subset\R^m$ is {\it basic} if for each continuous function 
$f:K\to\R$ there exist continuous functions $g_1,\dots,g_m:\R\to\R$ such that
$f(x_1,\dots,x_m)=g_1(x_1)+\dots+g_m(x_m)$ for each point
$(x_1,\dots,x_m)\in K$.

\smallskip
{\bf Theorem~1.} [St89] {\it Any $n$-dimensional compactum is basically 
embeddable into $\R^{2n+1}$ and, for $n>1$, is not basically embeddable into 
$\R^{2n}$.}

\smallskip
It is interesting to compare this theorem with the N\"obeling-Menger-Pontryagin 
theorem on embeddability of any $n$-dimensional compact space into $\R^{2n+1}$
and the example of an $n$-dimensional polyhedron non-embeddable into $\R^{2n}$. 

Obviously, $K$ is basically embeddable into $\R$ if and only if $K$
is topologically embeddable into $\R$. 
It follows from Theorem~1 that a compactum $K$ is basically embeddable into 
$\R^m$ for $m>2$ if and only if $\dim K<m/2$. 
Thus, the only remaining case is $m=2$ (Sternfeld's problem).

%\bigskip
%\centerline{\uppercase{\bf Basic embeddings into the plane}}

\bigskip
{\bf Basic embeddings into the plane}

%The problem of characterization of compact basic {\it subsets} of the plane 
%was stated back in [Arn58] and was solved in [Ste89]: a compact subset $K$ 
%of the plane is basic if and only if $E^n(K)=\emptyset$ for some $n$. Here
%$$E(Z)=\{z\in Z: \ |Z\cap p_x^{-1}p_xz|\ge2
%\text{ and }|Z\cap p_y^{-1}p_yz|\ge2\} $$
%and $p_x$, $p_y$ are the projections on the coordinate axes in the plane. 

A subset $K$ of $\R^2$ is {\it basic} if for each continuous 
function $f\colon K\to\R$ there exist continuous functions $g,h\colon\R\to\R$ 
such that $f(x,y)=g(x)+h(y)$ for each point $(x,y)\in K$.

Let us present the characterization of arcwise connected compacta basically 
embeddable into the plane [Sk95] (this is a partial solution of Sternfeld's 
problem). 
We formulate the criterion first for graphs and then for the general case.
A conjecture on embeddability of (not necessarily arcwise connected) 
connected compacta into the plane can be found in~[Sk95].
Compacta used in the statements are defined after the statements.  

\smallskip
{\bf Theorem~2.}~[Sk95] 
{\it A finite graph $K$ is basically embeddable into the plane if and only if
any of the following two equivalent conditions holds:

(a) $K$ does not contain subgraphs homeomorphic to $S, C_1, C_2$ (fig.~1), that 
is, a circle, a five-point star, and a cross with branched endpoints;

(b) $K$ is contained in one of the graphs $R_n$, $n=1,2,3,\dots$ (fig.~2).} 

%(cf. [Cla34, Cla37])

\begin{figure}
%[b]
  \begin{center}
    \begin{tabular}{c@{\qquad}c@{\qquad}c}
      \mpfile{gr}{1}
      &
      \mpfile{gr}{2}
      &
      \mpfile{gr}{3}
      \\[3mm]
      $S^1$& $T_5$ & $C$
    \end{tabular}
  \end{center}
  \caption{}\label{pic2}
\end{figure}

\begin{figure}
  \begin{center}
    \begin{tabular}{c@{\qquad}c@{\qquad}c}
      \raisebox{12.9mm}{\mpfile{gr}{31}}
      &
      \raisebox{3.87mm}{\mpfile{gr}{32}}
      &
      \mpfile{gr}{33}
      \\[3mm]
      $R_1$& $R_2$ & $R_3$\\[5mm]
      \raisebox{12.9mm}{\mpfile{gr}{34}}
      &
      \raisebox{3.87mm}{\mpfile{gr}{35}}
      &
      \mpfile{gr}{36}
      \\[3mm]
      $F_1$& $F_2$ & $F_3$
    \end{tabular}
  \end{center}
  \caption{}\label{pic3}
\end{figure}

Let $F_1$~be a triod. 
The graph $F_{n+1}$ is obtained from $F_n$ by branching 
%each of 
its endpoints (fig.~2). 
The graph $R_n$ is obtained from $F_n$ by by adding a hanging edge to each 
non-hanging vertex. 

\smallskip
{\bf Theorem~3.} [Sk95]
{\it An arcwise-connected compactum $K$ is basically embeddable into the plane 
if and only if it is locally connected (i.e., is a Peano continuum) and any of 
the two following (equivalent) conditions hold:

(1) $K$ does not contain $S^1, C_2, C_4, B$ as subcompacta and contains only 
finitely many subcontinua $F_n, H_n$ (fig.~1,2,3);

(2) $K$ does not contain any of the continua
$S^1, C_1, C_2, C_3, B, F$, $H_+, H_-, h_+, h_-$ (fig.~1,3,4).} 

% (3) there exist a filtration $I_1\subset\dots\subset I_n=K$ and a
%positive integer $M$ such that for any arc $s$ of order $l-1$
%($s=I_1$ if $l=1$), denoting by $\{s_m\}$ all the arcs of order at
%least $l$ that intersect $s$ and setting $R=s\cap(\bigcup_m s_m)$, we have:
%(I) $R\subset\Int s$;
%(C) at most two of the arcs $s_m$ are attached to each point of $s$,
%and if there are exactly two, then one of them is hanging (i.e., none
%of the arcs of order $l+1$ or greater intersects it), and
%the point to which these two arcs are attached has a neighborhood in $K$
%homeomorphic to the cross;
%(B) $\Cl R$ is nowhere dense in $S$;
%(H) $h^MQ=\emptyset$, where
%$Q=\{x\in R\ |\ \text{if }x\in s_m,\text{ then }s_m\text{ is not hanging}\}$
%and $hZ$ is a set obtained from $Z$ by removing points isolated in $Z$.

\begin{figure}
  \begin{center}
    \begin{tabular}{c@{\qquad}c@{\qquad}c}
      \mpfile{gr}{41}
      &
      \mpfile{gr}{42}
      &
      \mpfile{gr}{43}
      \\[3mm]
      $C_4$& $B$ & $H_2$
    \end{tabular}
  \end{center}
  \caption{}\label{pic4}
\end{figure}

\begin{figure}
  \begin{center}
    \begin{tabular}{c@{\qquad}c}
      \mpfile{gr}{51}
      &
      \mpfile{gr}{52}
      \\[2mm]
      $C_3$& $F$
      \\[3mm]
      \mpfile{gr}{53}
      & 
      \mpfile{gr}{54}
      \\[2mm]
      $H_-$&$H_+$
      \\[3mm]
      \mpfile{gr}{55}
       &
      \mpfile{gr}{56}
      \\[2mm]
      $h_+$& $h_-$  
    \end{tabular}
  \end{center}
  \caption{}\label{pic5}
\end{figure}

\smallskip
%Let us introduce notation and definitions needed to formulate the 
%characterization of arcwise connected continua basically embeddable into the 
%plane. 
Let $I=[0;1]$.
A sequence of sets is called a {\it null-sequence} if their diameters tend to 
zero.
Define 

$\bullet$ $H_n$ to be the union of $I$ with a null-sequence of triods having 
endpoints attached to $I$ at points $3^{-l_1}+\dots+3^{-l_s}$, where 
$s\le n$ and $0<l_1<\dots<l_s$ are integers;  

$\bullet$ $C_3$ to be a cross with a null-sequence of arcs attached to
one of its branches and converging to its center; 

$\bullet$ $C_4$ to be a cross with a sequence of points converging to its 
center; 

$\bullet$ $B$ to be the union of the arc $I$ and a null-sequence of arcs 
attached to $(0;1)$ by their endpoints at rational points;  

$\bullet$ $F$ to be the union of $I$ with a null-sequence of sets $F_n$ each
having an endpoint attached to the point $1/n\in I$;

$\bullet$ $H_+$ ($H_-$) to be the union of $I$ with a null-sequence of
continua $H_n$ connected to the points $1/n\in I$ by arcs that intersect $H_n$ 
at the points $1\in I\subset H_n$ ($0\in I\subset H_{n-1}$, respectively);   

$\bullet$ $h_+$ ($h_-$) to be obtained from a null-sequence of continua $H_n$ 
by pasting together the points $1\in I\subset H_n$ and $0\in I\subset H_{n-1}$
($0\in I\subset H_n$ and $1\in I\subset H_{n-1}$, respectively).

%Denote by $\Int J$ an arc $J$ without its endpoints. 
%(obviously, the topological type of $B$ 
%does not depend on possible variations in the construction); 
%We say that an arc $s$ {\it connects} points $A$ and $B$ if $s\cap A$ and 
%$s\cap B$ are different endpoints of $s$. 
%A sequence $I_1\subset I_2\subset\dots\subset I_n=K$
%such that $I_{l+1}$ is obtained by attaching to $I_l$
%a null-sequence of arcs, each at one endpoint, is called
%a {\it filtration}. These arcs will be called the
%{\it arcs of order $l$}.

\smallskip
%\bigskip
%{\bf Basic embeddings into products of graphs}

An embedding $K\subset X\times Y$ is {\it basic} if for any continuous function 
$f:K\to\R$ there exist continuous functions $g:X\to\R$, $h:Y\to\R$ such that
$f(x,y)=g(x)+h(y)$ for any point $(x,y)\in K$. 

Denote by $T_n$ an $n$-od, i.e., an $n$-pointed star. 
A vertex of a graph $K$ is called {\it horrible} if its degree is greater
than~4 and {\it awful} if its degree is equal to~4 and
it is not an endpoint of a hanging edge. 
The {\it defect} of a graph $K$ is the sum
$\delta(K)=(degA_1-2)+\dots+(degA_k-2)$, where $A_1,\dots ,A_k$
are all the horrible and awful vertices of~$K$.

\smallskip
{\bf Theorem~4.} [Ku99] {\it A finite graph $K$ admits a basic embedding 
$K\subset\R\times T_n$ if and only if $K$ is a tree and either $\delta(K)<n$ or
$\delta(K)=n$ and $K$ has a horrible vertex with a hanging edge.}

\bigskip
\centerline{\uppercase{\bf Basic planar sets}}

\smallskip
The material is presented as a sequence of problems, which is peculiar not only 
to Zen monasteries but also to elite mathematical education (at least in 
Russia). 
Difficult problems are marked by a star, and unsolved problems by two stars.
If the statement of a problem is an assertion, then it is required to prove
this assertion.

\bigskip
{\bf Discontinuously basic subsets.}

\smallskip
{\bf 1.} (a) Is it true that for any four numbers $f_{11},f_{12},f_{21},f_{22}$
there exist four numbers $g_1,g_2,h_1,h_2$ such that $f_{ij}=g_i+h_j$ for each
$i,j=1,2$?

(b) Andrey Nikolaevich and Vladimir Igorevich play the 'Dare you to decompose!'
game.
Some cells of chessboard are marked.
A. N. writes numbers in the marked cells as he wishes.
V. I. looks at the written numbers and chooses (as he wishes) 16 numbers
$a_1,\dots,a_8,b_1,\dots,b_8$ as 'weights' of the columns and the lines.
If each number in a marked cell turns out to be equal to the sum of
weights of the line and the row (of the cell), then V. I. wins, and in the
opposite case (i.e., when the number in at least one marked cell is not
equal to the sum of weights of the line and the row) A. N. wins.

Prove that V. I. can win no matter how A. N. plays if and only if there does 
not exist a closed route of a rook starting and turning only at marked cells
(the route is not required to pass through each marked cell).

\smallskip
Let $\R^2$ be the plane with a fixed coordinate system.
Let $x(a)$ and $y(a)$ be the coordinates of a point $a\in\R^2$.
An ordered set (either finite or infinite) $\{a_1,\dots,a_n,\dots\}\subset\R^2$
is called an {\it array} if for each $i$ we have $a_i\neq a_{i+1}$ and
$x(a_i)=x(a_{i+1})$ for even $i$ and $y(a_i)=y(a_{i+1})$ for odd $i$.
It is not assumed that points of an array are distinct.
An array is called {\it closed} if $a_1=a_{2l+1}$.

\smallskip
{\bf 2.} Consider a closed array $\{a_1,\dots,a_n=a_1\}$.
A {\it decomposition} for such an array is an assignment of numbers at the
projections of the points of the array on the $x$-axis and on the $y$-axis.
Is it possible to put numbers $f_1,\dots,f_n\in\R$, where $f_1=f_n$, at the
points of the array so that for each decomposition there exists an $f_i$
that is not equal to the sum of the two numbers at $x(a_i)$ and $y(a_i)$?

\smallskip
A subset $K\subset\R^2$ is called {\it discontinuously basic} if for each
function $f:K\to\R$ there exist functions $g,h:\R\to\R$ such that
 $f(x,y)=g(x)+h(y)$ for each point $(x,y)\in K$.

\smallskip
{\bf 3.} (a) The segment $K=0\times[0;1]\subset\R^2$ is discontinuously basic.

(b) The cross $K=0\times[-1;1]\cup[-1;1]\times0\subset\R^2$
is discontinuously basic.

(c) {\it A criterion for a subset of the plane to be
discontinuously basic.}
A subset of the plane is discontinuously basic if and only if it does
not contain any closed arrays.

\smallskip
{\bf 4.**} Given a set of marked unit cubes in the cube $8\times8\times8$, 
how can we see who wins in the 3D analogue of the `Dare you to decompose!' 
game?
In this analogue V. I. tries to choose 24 numbers
$a_1,\dots,a_8,b_1,\dots,b_8,c_1,\dots,c_8$ so that the number at the unit cube
$(i,j,k)$ would be equal to the sum $a_i+b_j+c_k$ of the three weights.

\smallskip
{\bf 5.**} (a) Define discontinuous basic subsets of the 3-space.
Discover and prove the 3D analogue of the above criterion.

(b) The same for higher-dimensional case.

\bigskip
{\bf Solutions.}

\small

\smallskip
{\bf 1.} (a) It is not true. If $f_{ij}=g_i+h_j$ for each $i,j=1,2$, then
$f_{11}+f_{22}=f_{12}+f_{21}$, but this is false for some numbers $f_{ij}$.

(b) The statement `only if' follows from the problem 2.
Let us prove the `if' part by induction on the number of the marked cells.
If only one cell is marked then we are done. 
Let $K$ be the set of centres of the marked cells. 
The set $E(K)$ is defined in the following subsection after Problem 9. 
The set $K$ does not contain any closed array, therefore $\#E(K)<\#K$.
So by the induction hypothesis V. I. can win for $E(K)$. 
Each cell from $K-E(K)$ is the only marked cell on its line or column, 
thus V. I. can choose the remaining weights for $K$.

\smallskip
{\bf 2.} Yes, it is. If every $f_i$ is equal to the sum of two numbers
at $x(a_i)$ and $y(a_i)$, then $f_1-f_2+f_3- \dots -f_{n-1}=0$, 
but this is false for some numbers $f_i$.

\smallskip
{\bf 3.} (a)  Set $h(y)=f(0,y)$ and $g(x)=0$.

(b) Set $g(x)=f(x,0)$ and $h(y)=f(0,y)-f(0,0)$.

(c)  The statement `only if' follows from the problem 2.
Let us prove the `if' part. Consider a function $f:K\to \R$. 
Our aim is to construct functions $g$ and $h$ so that $f(x,y)=g(x)+h(y)$. 
Two points $a,b\in K$ are called {\it equivalent} if there is an array 
$\{a=a_1,\dots,a_n=b\}\subset K$. 
Now take an equivalence class $K_1\subset K$. 
Define function $g:x(K_1)\to\R$ and $h:y(K_1)\to\R$ in the following way. 
Take any point $a_1\in K_1$ and set $g(x(a_1))=f(a_1)$ and $h(y(a_1))=0$.
If $\{a_1,a_2,\dots,a_{2l}\}$ is an array, then set 
$$h(y(a_{2l})):=f(a_{2l})-f(a_{2l-1})+\dots -f(a_1)\quad\text{and}
\quad g(x(a_{2l})):=f(a_{2l-1})-f(a_{2l-2})+\dots+f(a_1).$$
If $\{a_1,a_2,\dots,a_{2l+1}\}$ is an array, then set 
$g(x(a_{2l+1})):=f(a_{2l+1})-f(a_{2l})+\dots+f(a_1)$
($h(y(a_{2l+1}))$ is already defined). 
Make this construction for each equivalence class.
Then set $g=0$ and $h=0$ at all other points of $\R$. 

\normalsize

\bigskip
{\bf Continuously basic subsets.}

A subset $K\subset\R^2$ is called {\it (continuously) basic} if for each
continuous function $f:K\to\R$ there exist continuous functions $g,h:\R\to\R$
such that $f(x,y)=g(x)+h(y)$ for each point $(x,y)\in K$.

\smallskip
{\bf The Arnold problem.} {\it Which subsets of the plane are basic?} [Ar58]

\smallskip
In order to approach a solution consider some examples. 

\smallskip
{\bf 6.} (a) A closed array is not basic.

(b) The segment $K=0\times[0;1]\subset\R^2$ is basic.

(c) The cross $K=0\times[-1;1]\cup[-1;1]\times0\subset\R^2$ is basic.

(d) The graph $V$ of the function $y=|x|$, $x\in[-1;1]$ is basic.

\smallskip
A sequence of points $\{a_1,\dots,a_n,\dots\}\subset\R^2$
{\it converges to a point $a\in\R^2$} if for each $\eps>0$ there exists an
 integer $N$ such that for each $i>N$ we have $|a_i,a|<\eps$.

\smallskip
{\bf 7.} 
(a) If a subset of the plane is basic, then it is discontinuously basic.

(b) A {\it completed array} is the union of a point $a_0\in\R^2$ with an
infinite array $\{a_1,\dots,a_n,\dots\}\subset\R^2$ of distinct points which
converges to the point $a_0$. 
Prove that any completed array is not basic.
(Note that it is discontinuously basic).

(c) Let $[a,b]$ be the rectilinear arc which connects points $a$ and $b$.
Prove that the cross $K=[(-1,-2),(1,2)]\cup[(-1,1),(1,-1)]$ is not basic.

(d) Let $m_{ij}=2-3\cdot2^{-i}+j\cdot2^{-2i}$.
Consider the set of points $(m_{i,2l},m_{i,2l})$ and
$(m_{i,2l},m_{i,2l-2})$, where $i$ varies from 1 to $\infty$ and
$l=1,2,3,\dots,2^{i-1}$.
Prove that this subset of the plane does not contain any infinite arrays but
contains arbitrary long arrays.

(e) The union of the set from the previous problem and the point
 $(2,2)$ is not basic.

\smallskip
A subset $K\subset\R^2$ of the plane is {\it closed}, if for each sequence 
$a_i\in K$ converging to a point $a$ this point belongs to $K$.

\smallskip
{\bf 8.} A subset $K\subset\R^2$ of the plane is closed if and only if for each
point $a\not\in K$  there exists $\eps>0$ such that if for a point $b$ of the
plane we have $|a,b|\le\eps$, then $b$ does not belong to $K$.
%Prove the equivalence of both definitions of the closed subset.

\smallskip
{\bf The Sternfeld criterion for being a basic subset.} {\it  A closed
bounded subset $K\subset\R^2$ of the plane is basic if and only if
$K$ does not contain arbitrary long arrays.}

\smallskip
{\bf 9.}
(a) The criterion is false without the assumption that $K$ closed.

(b) The criterion is false without the assumption that $K$ bounded.

(c)** Find a criterion of being a basic subset for closed (but unbounded) 
subsets of the plane.  

\smallskip
Suppose that $K$ is a subset of $\R^2$.
For every point $v\in K$ consider the pair of lines passing through $v$ and
parallel to the $x$-axis and the $y$-axis.
If one of these two lines intersects $K$ only at point $v$, we colour $v$ in 
white.
Define $E(K)$ as the set of noncoloured points of $K$:
$$E(K)=\{v\in K:\ |K\cap(x=x(v))|\ge2\text{ and }|K\cap(y=y(v))|\ge2\}.$$
Let $E^2(K)=E(E(K))$, $E^3(K)=E(E(E(K)))$ etc.

\smallskip
{\bf 10.} A subset $K$ of the plane does not contain arbitrary long arrays 
if and only if $E^n(K)=\emptyset$ for some $n$.

\smallskip
{\bf 12.} (a)* Give an elementary proof that if $K$ is a closed bounded subset 
of $\R^2$ and $E(K)=\emptyset$, then $K$ is basic [Mi09].

Hint. It can be proven that for piecewise-linear maps $f$
there is a decomposition $f(x,y)=g(x)+h(y)$  with $|g|+|h|<5|f|$. 

(b)* Prove the `if' part of the criterion without using the 
functional spaces 
%Banach Inverse Operator Theorem 
as below.

Hint. Same as above with $|g|+|h|<C_n|f|$,
where $C_n$ depends only on that $n$ for which $E^n(K)=\emptyset$.

\smallskip
{\bf 11.}  A subset $K\subset\R^3$ is called {\it (continuously) basic} if for 
each continuous function $f:K\to\R$ there exist continuous functions 
$g,h,l:\R\to\R$ such that $f(x,y,z)=g(x)+h(y)+l(z)$ for each point $(x,y,z)\in K$.

(a) The `hedgehog'
$0\times0\times[-1;1]\cup0\times[-1;1]\times0\cup
[-1;1]\times0\times0\subset\R^3$ is basic.

(b) The set of 4 points $(0,0,0)$; $(1,1,0)$; $(0,1,1)$; $(1,0,1)$ is  basic.
(But $E^n(K)\neq\emptyset$ for each $n$, see below.)

(c)* Define $E(K)$ analogously to the above, only instead of lines use planes 
orthogonal to the axes:
$$E(K)=\{v\in K:
\ |K\cap(x=x(v))|\ge2,\ |K\cap(y=y(v))|\ge2\text{ and }|K\cap(z=z(v))|\ge2\}.$$
Let $K$ be a closed bounded subset of $\R^3$. 
Prove that if $E^n(K)=\emptyset$ for some $n$, then $K$ is basic 
[St89, Lemma 23.ii].

\bigskip
\newpage
{\bf Solutions.}

\small
\smallskip
{\bf 6.}
(a) If an array $A=\{a_1,\dots,a_{2l+1}\}$ is basic, then
$f(a_1)-f(a_2)+\dots +f(a_{n-2})-f(a_{2l})=0$. 
But this is false for some functions $f$. Cf. problem 2.

(b),(c) Analogously to problems 3a,3b.

(d) Take $h(y)=0$ and $g(x)=f(x,y)$.

\smallskip
{\bf 7.} (a) If the subset is not discontinuously basic, then it contains a 
closed array. 
Hence the statement follows by extension of $f$ on the subset and using 
problem 6a. 

(b) Define function $f$ by $f(a_n)=\frac{(-1)^n}n$. Suppose that
 $f(x,y)=g(x)+h(y)$ for some $g$ and $h$. 
Then 
$$f(a_1)-f(a_2)+f(a_3)-f(a_4)+ \dots-f(a_{2l})=h(y(a_1))-h(y(a_{2l})).$$
Since $\lim_{l\to\infty}h(y_{2l})$ exists and equals to $h(y(a_0))$,
it follows that $\sum_{i=1}^{2l} (-1)^i f(a_i)$ converges
when $l\to\infty$, which is a contradiction.

(c) The cross contains a completed array 
$$a_{4k+1}=(-2^{-2k},2^{-2k}),\ a_{4k+2}=(2^{-2k-1},2^{-2k}),
\ a_{4k+3}=(2^{-2k-1},-2^{-2k-1}),\ a_{4k+4}=(-2^{-2k-2},-2^{-2k-1}).$$
Define a function $f$ on this array using problem 7.b and then extend it 
(e.g. piecewise linearly) to the cross. 
Then there are no functions $g$ and $h$ such that $f(x,y)=g(x)+h(y)$.

(d) For every $i$ the set 
 ${(m_{i,2l},m_{i,2l})}_{l=1}^{2^{i-1}}\cup 
{(m_{i,2l},m_{i,2l-2})}_{l=1}^{2^{i-1}}$
is an array of $2^i$ points.

(e) Define a function $f$ by 
$$f((m_{i,2l},m_{i,2l})):=2^{-i}\quad\text{and}
\quad f(m_{i,2l},m_{i,2l-2}):=-2^{-i}.$$ 
If $f(x,y)=g(x)+h(y)$ for some $g$ and $h$, then for every $i$ using array 
of points $(m_{i,2l},m_{i,2l})$ and $(m_{i,2l},m_{i,2l-2})$, 
where $l=1,2,3,\dots 2^{i-1}$, we obtain  
$h(2-\frac 3{2^i})-h(2-\frac 2{2^i})=1$.
 This contradicts to the continuity of $h$.

\smallskip
{\bf 8.} Let us prove the `only if' part.
Let $K$ be a closed  subset of the plane.
Suppose that for some point $a=(x,y)\not\in K$ and for each $\eps=\frac 1n>0$
there exists a point $a_n\in K$ (at least one) such that $|a,a_n| \le \frac 1n$. 
The sequence of points $a_n\in K$ converges to the point $a$, thus $a\in K$. 
 Contradiction.

Now let us prove the `if' part.
Suppose that a sequence $a_n$ converges to a point $a$ and the point
$a=(x,y)$ is not in $K$. There exists $\eps>0$ such that
 for every point $a_n \in K$ the distance $|a,a_n|>\eps$. 
This is a contradiction.

\smallskip
{\bf 9.} (a) Any infinite array $A$ not containing closed arrays and converging
to a point $a\not\in A$ is basic. 
This follows because each function defined on $A$ is continuous.  

(b) A counterexample is
$\{(k,k)\}_{k=1}^\infty\cup\{(k,k-1)\}_{k=1}^\infty$.

\smallskip
{\bf 10.} Let us prove the `only if' part. 
 Suppose that $E^n(K)\neq\emptyset$ for each $n$.
For each $n$ take a point $a_0\in E^n(K)$. Then there exist points $a_{-1},a_1\in
E^{n-1}(K)$ such that $x(a_{-1})=x(a_0)$ and $y(a_1)=y(a_0)$.
Analogously there  exist points $a_{-2},a_2\in E^{n-2}(K)$ such that
$\{a_{-2},a_{-1},a_0,a_1,a_2\}$ is an array. Analogously we construct
an array of $2n+1$ points in $K$, which is a contradiction. 
%Thus $E^n(K)=\emptyset$ for some $n$.

Let us prove the `if' part. 
Suppose that $K$ contains an array of $2n+1$ points
 $\{a_{-n},\dots,a_0,\dots,a_{n}\}$. Then there is an array of $2n-1$ points
 $\{a_{-n+1},\dots,a_{n-1}\}$ in $E(K)$. Analogously $a_0\in E^n(K)$. 
Thus if $E^n(K)=\emptyset$, then
 $K$ does not contain an array of $2n+1$ points.

\smallskip
{\bf 11.}
(a) For each functuion $f:K \to \R$ on $K$ define 
 $g(x):=f(x,0,0)$, $h(y):=f(0,y,0)-f(0,0,0)$ and $l(z):=f(0,0,z)-f(0,0,0)$. 

(b) Set $g(0)=f(0,0,0)$, $h(0)=0$, $l(0)=0$,
$$2g(1)=f(0,0,0)+f(1,1,0)+f(1,0,1)-f(0,1,1),\quad 
2h(1)=-f(0,0,0)+f(1,1,0)-f(1,0,1)+f(0,1,1)$$
$$\text{and}\quad
2l(1)=-f(0,0,0)-f(1,1,0)+f(1,0,1)+f(0,1,1).$$

\normalsize

\smallskip
{\bf Proof of the criterion for being a basic subset.}

\smallskip
Let $K$ be a closed bounded subset of the plane.
It is known that each continuous function $f:K\to \R$ is bounded.
A function $f:K\in\R$ is called {\it bounded}, if there exists a number $M$ 
such that $|f(x)|<M$ for every $x\in K$. 
For a bounded function $G:K\to\R$ denote $|G|:=\sup_{x\in K}|G(x)|$.

\smallskip
{\it Beginning of the proof of the `only if' part of the criterion.}
Assume to the contrary that $K$ contains arbitrary long arrays and is basic.
Choosing subsequences we may assume that points of each array are distinct.
Therefore for each $n$ there is an array $\{a^n_1,\dots,a^n_{2n+5}\}$ of $2n+5$
distinct points in $K$.

Then there exists continuous function 
$$f_n:K\to\R\quad\mbox{such that}\quad f_n(a^n_i)=(-1)^i\quad\mbox{and}\quad 
|f_n(x)|\le1\quad\mbox{for each}\quad x\in K.$$
(Indeed, first define such a continuous function $f:\R^2\to \R$. 
Denote $s=\min_{i<j}|a_i,a_j|$. 
Take $n$ disks with centers $a_i$ and radii $\frac s3$. 
Outside of these disks set $f=0$.
Inside the $i$-th disk take $f$ to be $(-1)^i$ in the center $a_i$, $0$ on the 
boundary and extend it linearly in the distance to $a_i$.  
Then restrict $f$ to $K\subset \R^2$.)

%Let $f:K\to\R$ and $g,h:\R\to\R$ be functions such that $|f-f_n|<1/2n$ and
%$f(x,y)=g(x)+h(y)$ for each $(x,y)\in K$. Then  by the problem 5(d) $|g|>n$.

Define integers $s_n$ and functions $F_n:K\to\R$ inductively as follows.
Set $s_0=1$ and $F_0=0$.
Suppose now that $F_{n-1}$ and $s_{n-1}$ are defined.
If $F_{n-1}$ is not representable as $G_{n-1}(x)+H_{n-1}(y)$, then we are done.
If it is representable in this way, then take 
$$s_n>s_{n-1}!(|G_{n-1}|+n)
\quad \text{ and } \quad  F_n=F_{n-1}+\dfrac{f_{s_n}}{s_{n-1}!}$$
It remains to prove that if we can construct in this way an infinite number of
$s_n$ and $F_n$, then the function
$$F=\lim\limits_{n\to\infty}F_n=\sum\limits_{n=1}^\infty\frac{f_{s_n}}{s_{n-1}!}$$
is not representable as $G(x)+H(y)$.

Assume to the contrary that $F(x,y)=G(x)+H(y)$ for some $G$ and $H$. 
It suffices to prove that $|G|>n$ for each $n$. 
For this it suffices to prove that $s_{n-1}!|G-G_{n-1}|>s_n$: then we would 
have $$|G|+|G_{n-1}|\ge|G-G_{n-1}|>\frac{s_n}{s_{n-1}!}>
|G_{n-1}|+n.$$

%\smallskip
{\bf Lemma.} {\it Let $m\ge 4$, 

$\bullet$ $K=\{a_1,\dots,a_{2m+5}\}$ be an array of $2m+5$ distinct points, 

$\bullet$ $f(a_1),\dots,f(a_{2m+5})$ numbers such that $|(-1)^i-f(a_i)|\le 1/m$, 

$\bullet$ $g(x(a_i)),h(y(a_i))$, $i=1,\dots,2m+5$, numbers such that 
$f_i=g(x(a_i))+h(y(a_i))$ for each $i$.

Then $\max_i |g(x(a_i))|>n$. }

\smallskip
{\it Proof.} We may assume that $a_1a_2\| Ox$. 
Then 
$$|(f(a_1)-f(a_2)+f(a_3)-f(a_4)+\dots-f(a_{2m+4})-(2m+4)|\le\frac{2m+4}m\le3.$$
Therefore $g(x(a_1))-g(x(a_{2m+4}))\ge(2m+4)-3>2m$. 
This implies the required inequality.
\qed

\smallskip
{\it Completion of the proof of the `only if' part of the criterion.}
We have
$$F-F_n=F-F_{n-1}-\frac{f_{s_n}}{s_{n-1}!}=
\frac{s_{n-1}!(F-F_{n-1})-f_{s_n}}{s_{n-1}!}.$$
Apply the Lemma to 
$$m=s_n,\quad a_i=a_i^{s_n},\quad f=s_{n-1}!(F-F_{n-1}),
\quad g=s_{n-1}!(G(x)-G_{n-1}(x)), \quad h=s_{n-1}!(H(y)-H_{n-1}(y)).$$
This is possible because $f(x,y)=g(x)+h(y)$ and (since $s_n-1>s_{n-1}$ for 
$n>2$)  
$$|f-f_{s_n}|=s_{n-1}!|F-F_n|<\frac1{(s_n-1)\cdot s_n}\sum\limits_{k=0}^\infty
\frac1{(s_n+1)\cdot\dots\cdot s_{n+k}}<
\frac1{(s_n-1)\cdot s_n}\sum\limits_{k=0}^\infty\frac1{2^k}<\frac1{s_n}.$$
By Lemma we obtain $s_{n-1}!|G-G_{n-1}|>s_n$.
\qed

\smallskip
{\it Proof of the criterion. 
\footnote{This proof is not elementary, is not used in the sequel and could 
be omitted.}
}
The proof is based on a reformulation of the property of being a basic subset 
in terms of {\it bounded linear operators} in {\it Banach functional spaces}.  
Denote by $C(X)$ the space of continuous functions on
$X$ with the norm $|f|=\sup\limits\{|f(x)|\ :\ x\in X\}$.  
In this proof denote by $pr_x(a)$ and $pr_y(a)$ the projections of a point 
$a\in K$ on the coordinate axes.

For $K\subset I^2:=[0;1]\times[0;1]$ define a map ({\it linear superposition 
operator}) 
$$\phi\colon C(I)\oplus C(I)\to C(K)\quad\text{by}\quad
\phi(g,h)(x,y):=g(x)+h(y).$$
Clearly, the subset $K\subset I^2$ is basic if and only if $\phi$ is 
surjective, or equivalently, epimorphic. 

Denote by $C^*(X)$ the space of {\it bounded linear functions} $C(X)\to\R$
with the norm $|\mu|=\sup\{|\mu(f)|\ :\ f\in C(X),\ |f|=1\}$.
For a subset $K\subset I^2$ define a map ({\it dual linear superposition 
operator}) 
$$\phi^*\colon C^*(K)\to C^*(I)\oplus C^*(I)\quad\text{by}\quad
\phi^*\mu(g,h):=(\mu(g\circ pr_x),\mu(h\circ pr_y)).$$
Since $|\phi^*\mu|\le2|\mu|$, it follows that $\phi^*$ is bounded.
By duality, $\phi$ is epimorphic if and only if $\phi^*$ is monomorphic.
\footnote{We remark that  $\phi^*$ can be injective but not monomorphic.
In other words not only some linear relation on $\im\phi$ can force it to
be strictly less than $C(K)$.
\newline
If an embedding $K\subset \R^2$ is basic, then we can prove that $\phi^*$ is 
monomorphic without use of $\phi$ as follows.  
Define a linear operator
$$\Psi\colon C^*(I)\oplus C^*(I)\to C^*(K)\quad\text{by}\quad
\Psi(\mu_x,\mu_y)(f)=\mu_x(g)+\mu_y(h),$$
where $g,h\in C(I)$ are such that $g(0)=0$ and $f(x,y)=g(x)+h(y)$ for
$(x,y)\in K$.  
Clearly,
$\Psi\phi^*=\id$ and $\Psi$ is bounded, hence $\phi^*$ is monomorphic.
} 

It is clear that $\phi^*$ is monomorphic if and only if

{\it (*) there exists $\varepsilon>0$ such that
$|\phi^*\mu|>\varepsilon|\mu|$ for each unzero $\mu\in C^*(K)$.}

%(because this condition ensures that $\im\phi^*$ is closed).
%Thus $K$ is basic if and only if the property $(\Phi)$ holds.

We leave as an excercise the proof that (*) implies the abcense of arbitrarily 
large arrows. 
(This proves the `only if' part of the criterion, for which we already have an 
elementary proof.)  

So it remains to prove that $E^n(K)=\emptyset$ implies the condition (*).
We present the proof for $n\in\{1,2\}$. 
The proof for arbitrary $n$ is analogous.
We use the following non-trivial fact: 
{\it $C^*(X)$ is the space of $\sigma$-additive 
regular real valued Borel measures on $X$}  
(in the sequel we call them simply `measures').
We have 
$$\phi^*\mu=(\mu_x,\mu_y),\quad\text{where}\quad\mu_x(U)=\mu(pr_x^{-1}U)
\quad\text{and}
\quad \mu_y(U)=\mu(pr_y^{-1}U)\quad\text{for each Borel set}\quad U\subset I.$$
If $\mu=\mu^+-\mu^-$ is the decomposition of a measure $\mu$ into its
positive and negative parts, then $|\mu|=\bar\mu(X)$,
where $\bar\mu=\mu^++\mu^-$ is the absolute value of $\mu$.

Let $D_x$ ($D_y$) be the set of points of $K$ which are
not shadowed by some other point of $K$ in $x$- ($y$-) direction.
Take any measure $\mu$ on $K$ of the norm 1.

If $n=1$, then 
$$E(K)=\emptyset,\quad\text{then}\quad D_x\cup D_y=K,\quad\text{so}\quad 
1=\bar\mu(K)\le\bar\mu(D_x)+\bar\mu(D_y).$$
Therefore without loss of generality, $\bar\mu(D_x)\ge1/2$.
Since the projection onto the $x$-axis is injective over $D_x$, it follows 
that $|\mu_x|\ge1/2$, thus the required assertion holds for $\varepsilon=\frac 12$.

If $n=2$, then 
$$E(E(K))=\emptyset,\quad\text{hence}\quad D_x\cup D_y=K-E(K),\quad\text{so}
\quad E(D_x\cup D_y)=\emptyset.$$
\quad 
In the case when $\bar\mu(E(K))<3/4$ we have
$\bar\mu(D_x\cup D_y)>1/4$ and without loss of generality $\bar\mu(D_x)>1/8$.
Then as for $n=1$ we have $|\mu_x|>1/8$, thus (*) holds for $\varepsilon
=\frac 18$.

In the case when $\bar\mu(E(K))\ge3/4$ we have $\bar\mu(K-E(K))\le1/4$.
By the case $n=1$ above without loss of generality 
$\bar\mu_x(pr_x(E(K)))\ge\bar\mu(E(K))/2$.
Hence $|\mu_x|\ge\frac12\cdot\frac34-\frac14=\frac18$, thus (*) holds for 
$\varepsilon=\frac 18$.
\qed

\bigskip
{\bf Smoothly basic subsets of the plane.}

Let $K$ be a subset of the plane $\R^2$. A function $f:K\to\R$ is called 
{\it differentiable} if for each point $z_0\in K$
there exist a vector $a\in\R^2$ and infinitesimal function
$\alpha:\R^2\to\R$ such that for each point $z\in K$
$$f(z)=f(z_0)+a\cdot(z-z_0)+\alpha(z-z_0)|z,z_0|.$$
\quad 
Here the dot denotes scalar product of vectors 
$a=:(f_x,f_y)$ and $z-z_0=:(x,y)$,
i.e. $a\cdot(z-z_0)=xf_x + yf_y$. 
A function $\alpha:\R^2\to\R$ is {\it infinitesimal}, if for each number 
$\eps >0$ there exists a number $\delta>0$ such that for each point
$(x,y)\in \R^2$ 
$$\text{if}\quad \sqrt{x^2+y^2}<\delta,\quad\text{then}
\quad|\alpha(x,y)|<\eps.$$
\quad
Let $V$ be the graph of the function $y=|x|$, where $x\in[-1;1]$. 
A function $f:V\to\R$ is differentiable if and only if $f(x,|x|)$ is 
differentiable on the segments $[-1;0]$ and $[0;1]$.

A subset $K\subset\R^2$ of the plane is called {\it differentiably basic} if 
for each differentiable function $f:K\to\R$ there exist differentiable 
functions $g:\R\to\R$ and $h:\R\to\R$ such that $f(x,y)=g(x)+h(y)$ for each 
point $(x,y)\in K$.

\smallskip
{\bf 13.} (a) (b) (c) Solve the analogues of problem 6 for differentiably
 basic sets.

\smallskip
{\bf 14.}
(a) The graph $V$ is differentiably basic.

(b) $W:=(V-(2,0))\cup(V+(2,0))$ is not differentiably basic. 

(c) The broken line whose consecutive vertices are $(-2,0)$, $(-1,1)$, $(0,0)$, $(1,1)$
and $(2,0)$ is not differentiably basic.
(Note that it is continuously basic).

(d) The completed array
$\{([\frac{n+1}2]^{-1/2},[\frac n2]^{-1/2})\}_{n=2}^{\infty}\cup\{(0,0)\}$
is not differentiably basic. (Note that it is also not continuously basic.)

(e) The completed array
 $\{(2^{-[\frac{n+1}2]},2^{-[\frac n2]})\}_{n=1}^{\infty}\cup\{(0,0)\}$
is differentiably basic. (Note that it is not continuously basic.)

(f) (I. Shnurnikov) 
The cross $K=[(-1,-2),(1,2)]\cup[(-1,1),(1,-1)]$ is not differentiably 
basic.
(This assertion and Conjecture 15a imply that the property of being 
differentably basic is not hereditary.)

(g) If a graph is basically embeddable in the plane, then it is differentiably 
basically embeddable in the plane.  
(This is non-trivial because the plane contains graphs which are 
basic but not differentaibly basic and vice versa.) [RZ06]

\smallskip
{\bf 15.**} {\bf Conjectures.} 
(a) (I. Shnurnikov) A completed array $\{a_n\}_{n=1}^\infty\cup\{(0,0)\}$ 
is differentiably basic if and only if the sequence 
$\frac{\sum\limits_{n=k}^\infty|a_n|}{|a_k|}$ is bounded. 

(b) The subset 
$\{(t^2,\frac {t^2}{(1+t)^2})\}_{t\in[-\frac 12;\frac 12] }$
of the plane is not differentiably basic. 

Hint. One can try to prove this analogously to 14f. Cf. [Vo81, Vo82]. 

(c) A piecewise-linear graph in $\R^2$ is differentiably basic if and only if 
it does not contain arbitrary long arrays and for each two singular points $a$ 
and $b$ we have $x(a)\ne x(b)$ and $y(a)\ne y(b)$.
A point $a\in K$ is {\it singular} if the intersection of $K$ with each disk 
centered at $a$ is not a rectilinear arc.

\smallskip
It would be interesting to find a criterion of being differentiably basic for 
closed bounded subsets of the plane.
Apparently a simple-to-state criterion (analogous to the Sternfeld criterion) 
does not exist. 
Another interesting question: is there a continuous map $[0;1]\to\R^2$ whose 
image is differentiably basic but not basic?

\smallskip
{\bf 16.} Let $r\ge0$ be an integer and $K\in \R^2$ a subset. 
A function $f:K\to\R$ is called {\it $r$ times differentiable} if for each 
point $z_0\in K$ there exist a polynomial $\overline f(z)=\overline f(x,y)$ 
of degree at most $r$ of 2 variables $x$ and $y$
and an infinitesimal function $\alpha:\R^2\to\R$ such that
$f(z)=\overline f(z-z_0)+\alpha(z-z_0)|z,z_0|^r$ for each point $z\in K$. 
(This definition differs from the one generally accepted.)

(a) Functions differentiable zero times are exactly continuous functions,
and functions differentiable one time are exactly differentiable functions.

(b) For each positive integer $r$ define the property of being an $r$ times
 differentiably basic subset of the plane $\R^2$.

(c) For each integer $k\ge0$ there is a subset of the plane which is $r$ times
 differentiably basic for  $ r=0,1\dots k$  but is not $r$  times
differentiably basic for each $r>k$.

(d)** Find a criterion for graphs in $\R^2$ to be $r$ times differentiably
basic.

\bigskip
{\bf Solutions.}

\small

\smallskip
{\bf 13.}
(a), (b), (c) Analogously to problems 6(a), 3(a) and 3(b).

\smallskip
{\bf 14.}
(a) Take a differentiable function $f:V\to\R$.
Since $f$ is differentiable at $(0,0)$, it follows that
 there exist $a,b\in\R$ such that 
$$f(x,|x|)=f(0,0)+ax+b|x|+\alpha(x),\quad\text{where}
\quad \alpha(x)=o(\sqrt{x^2+|x|^2})\quad\text{when}\quad x\to0.$$
Take $h(y):=by$ and $g(x):=f(0,0)+ax+\alpha(x)$. 
Clearly, $h$ is  differentiable and $g$ is differentiable outside 0. 
Since $\alpha(x)=o(x)$ when $x\to0$, it follows that $g$ is differentiable
also at 0.  

(b) See 16c for $k=0$. 

(c) Suppose the broken line is differentiably basic.
The function $f(x,y)=xy$ is differentiable. 
We have $f(x,y)=g(x)+h(y)$, where both $g$ and $h$ are differentiable.
Then 
$$2-2d=f(1+d,1-d)+f(1-d,1-d)=g(1+d)+g(1-d)+2h(1-d)=
2g(1)+2h(1)-2h'(1)d+o(d).$$
Hence $h'(1)=1$. 
Analogously  
$$2d-2=f(-1+d,1-d)+f(-1-d,1-d)=g(-1+d)+g(-1-d)+2h(1-d)=
2g(-1)+2h(1)-2h'(1)d+o(d).$$
Hence $h'(1)=-1$. 
A contradiction. 

(d) Suppose that this completed array is differentiably basic.
Set $a_n=([\frac{n+1}2]^{-1/2},[\frac n2]^{-1/2}),$
 $f(a_n):=\frac{(-1)^n}n$, 
$n=2,3,\dots$.
 If $f(x,y)=g(x)+h(y)$ for some functions $g(x)$ and $h(y)$,
then the series 
$f(a_2)-f(a_3)+f(a_4)-\dots$ converges to $g(1)-g(0)$ 
(analogously to Problem 7b).  
This is a contradiction because the series 
$\frac 12 +\frac 13+ \frac 14 + \dots $ diverges.

(e) Without loss of generality assume that $f(0,0)=0$, then take 
$g(0)=0$ and $h(0)=0$.
Set 
$$h(2^{-k})=f(2^{-(k+1)},2^{-k})-f(2^{-(k+1)},2^{-(k+1)})+f(2^{-(k+2)},2^{-(k+1)})-
\dots, $$
$$g(2^{-k})=f(2^{-k},2^{-k})-f(2^{-(k+1)},2^{-k})+
f(2^{-(k+1)},2^{-(k+1)})-\dots,$$
where the right-hand sides are sums of alternating series. Now $g(x)$ and
$h(y)$ may be extended to differentiable functions $\R\to\R$.

(f) 
Define 
$$w(0)=w(4^{-i}+4^{-3i})=w(4^{-i})=0\quad\text{and}
\quad w(4^{-i}+4^{-3i-1})=2^{3i}\quad\text{for}\quad i=1,2,3,\dots.$$ 
Extend piecewice-linearly to obtain a function $w:[0;1]\to \R$. 
For $x\in[0;1]$ define $W(x)$ as the area under the graph of $w$ on $[0;x]$.  
(This is well-defined because this area is finite.)
Define $f(x,-x)=W(x)$ for $x\in[0;1]$ and $f(x,y)=0$ on the rest of the cross.

Clearly, $f$ is differentiable outside $(0,0)$. 
It is easy to check that $f$ is differentiable at $(0,0)$. 

Suppose that $f(x,y)=g(x)+h(y)$ for some differentiable functions $g$ and $h$.
Without loss of generality we assume that $g(0)=h(0)=0$.
The function $g$ is not differentiable at $x=1/4$ because for $0<d<\frac 14$ 
we have 
$$g\left(\frac14+d\right)-g\left(\frac14\right)=
W\left(\frac14+d\right)-W\left(\frac14\right)+
W\left(\frac1{4^2}+\frac d4\right)-W\left(\frac1{4^2}\right)+\dots >$$
$$>W\left(\frac1{4^{k+1}}+\frac d{4^k}\right)-W\left(\frac1{4^{k+1}}\right)=
\frac{2^{3k}\cdot 4^{-3k}}2\ge \frac {(4d)^{3/4}}2.$$
Here 

$\bullet$ the first equality is proved using two infinite arrays starting 
at points $(\frac14+d,-\frac14-d)$ and $(\frac14,-\frac14)$ and converging to 
the point $(0,0)$; 

$\bullet$ $k\ge0$ is such that $4^{-2k}\ge 4d>4^{-2(k+1)}$;  

$\bullet$ the first inequality follows because $W$ is a non-decreasing function; 

$\bullet$ the second inequality follows because 
$\frac d{4^k}>\frac1{4^{3(k+1)}}$; 

$\bullet$ the second equality follows by definition of $k$.

(In the same way one can prove that $g$ is not differentiable at $x=4^{-i}$ 
for each $i$.) 

\smallskip
{\bf 15.} 
(a) Hints. 
For the `only if' part use the idea of Problem 7b and prove that 
if $\sum\limits_{n=1}^\infty|a_n|=\infty$, then there is a sequence $b_n\to0$ 
such that $\sum\limits_{n=1}^\infty|a_n|b_n=\infty$. 

For the `if' part we may assume that numbers $x(a_i)$ are distinct, numbers 
$y(a_i)$ are distinct, $x(a_{2i})=x(a_{2i+1})$, $y(a_{2i})=y(a_{2i-1})$. 
If $f(0,0)=0$, define 
$$g(x(a_{2i})):=f(a_1)-f(a_2)+f(a_3)-\dots+f(a_{2i+1}),
\quad g(0):=\sum_{i=1}^{\infty}(-1)^if(a_i),$$  
$$h(y(a_{2i})):=-f(a_1)+f(a_2)-f(a_3)-\dots+f(a_{2i-2})\quad\text{and}
\quad h(0):=\sum_{i=1}^{\infty}(-1)^if(a_i).$$  
Prove that $g$ and $h$ are differentiable at 0.

\comment

Since $f$ is differentiable at points $(-1,1)$ and $(1,1)$, the 
following relations hold for sufficiently small $d>0$:
$$f(-1+d,1-d)-f(-1,1)=f_1d-f_2d+\alpha_{(-1,1)}(d,-d)|(d,-d)|,$$ 
$$ f(-1-d,1-d)-f(-1,1)=-f_1d-f_2d+\alpha_{(-1,1)}(-d,-d)|(-d,-d)|,$$
$$f(1+d,1-d)-f(1,1)=f_3d-f_4d+\alpha_{(1,1)}(d,-d)|(d,-d)|\quad\text{and}$$
$$f(1-d,1-d)-f(1,1)=-f_3d-f_4d+\alpha_{(1,1)}(-d,-d)|(-d,-d)|.$$
Also we have $f(x,y)=g(x)+h(y)$ and both $g(x)$, $h(y)$ are differentiable. 
Hence
$$f(-1+d,1-d)-f(-1,1)=g(-1+d)-g(-1)+h(1-d)-h(1)=
  g'(-1)d-h'(1)d+\alpha(d)d\quad\text{and}$$
$$ f(-1-d,1-d)-f(-1,1)=g(-1-d)-g(-1)+h(1-d)-h(1)=
 -g'(-1)d-h'(1)d+\alpha(d)d.$$
Therefore $h'(1)=f_2$ (and $g'(-1)=f_1$). 
Analogously $h'(1)=f_4$.
Thus $h'(1)=f_2=f_4$. 
But  for function $f(x,y)=xy$ we have $f_4=1\neq f_2=-1$.

\smallskip
To 14e, 15a:

Let us introduce some notation for arrows converging to $(0,0)$.
Denote
$$a_1 = (x_1,y_1), \quad a_2=(x_2,y_1),\quad a_3=(x_2,y_2), \dots$$
$$a_{2k}=(x_{k+1},y_k),\quad a_{2k+1}=(x_{k+1},y_{k+1}),\dots$$
$$f_n:=f (a_n),\quad g_k:=g(x_k),\quad h_k:=h(y_k),$$
$$f(0,0):=f_\infty, \quad g(0):=g_\infty, \quad h(0):=h_\infty$$

In our case it suffices to consider differentiability only at the point
$(0,0)$, because at other points any functions
$g(x)$ and $h(y)$ will be differentiable.

The function $f_n$ is differentiable at point $(0,0)$ if there exists
$a,b\in R$ for which

$f_{2k-1}=f_\infty+\frac{a}{2^k}+\frac{b}{2^k}+\frac{\alpha_k}{2^k}$ and 
$f_{2k}=f_\infty+\frac{a}{2^{k+1}}+\frac{b}{2^k}+\frac{\beta_k}{2^k}$

where $\alpha_k$ and $\beta_k$ are infinitesimal.

In order to prove that $M_a$ is $D$-basic we only need to consider
the case when $f_{\infty}=a=b=0$.
Indeed, denote $s(x,y)=f(x,y)-f_\infty-ax-by$, then
from an expansion of the function $s$ the expansion of the function $f$
is easily constructed.
Thus it suffices to prove the existence of the expansion $f(x,y)=g(x)+h(y)$
for a function $f$ such that
$f_{2k-1}=\frac {\alpha_k}{2^k}$ and $f_{2k}=\frac {\beta_k}{2^k}$,  
i. e. such that $f_k=\frac{\gamma_k}{(\sqrt2)^k}$.

Take $h_1=0$.
Sum up the equalities $(*)$ from the Proof of the Lemma on arrows
with alternative sign, we have

$g_k=f_1-f_2+f_3-... +f _ {2k-1}$.
Thus
$$g_{\infty}=-\sum_{i=1}^{\infty}(-1)^{i}f_i=
-\sum_{i=1}^{\infty}(-1)^{i}\frac{\gamma_i}{(\sqrt2)^i},$$

and this sum, obviously, converges.
Analogously

$$h_k=-f_1+f_2-f_3+...+f_{2n-2}
\quad\text {and}
\quad h_{\infty}=\sum_{i=1}^{\infty}(-1)^{i}f_i $$

Now it suffices to prove
differentiability of the constructed functions $g_k$ and $h_k$.
We have
$$|\frac{h_k-h_\infty}{1/2^k}|=|\sum_{i=2k-1}^{\infty}2^k(-1)^{i}f_i|
\le\sum_{i=2k-1}^{\infty}
|\frac {\gamma_i}{(\sqrt {2})^{i- (2k-1)}} | \le $$
$$\le
\sqrt2 \sum _ {i=1} ^ {\infty} \frac {\varepsilon_k} {(\sqrt {2}) ^ {i}} =
\frac {2} {\sqrt2-1} \varepsilon_k, $$

where $\varepsilon_k=max_{i\ge 2k-1}\{|\gamma_i|\}$.

Therefore
 $$h'(0)=|\lim\limits_{k\to\infty}\frac{h_k-h_\infty}{1/2^k}|
 \le\lim\limits_{k\to\infty}\frac{2}{\sqrt2-1}\varepsilon_k=0.$$

The differentiability of the function $g_k $ is proved analogously.

\endcomment

\smallskip
{\bf 16.} (a) It is clear. 

(b) A subset $K\subset\R^2$ is called {\it $r$ times differentiably basic} if
for each $r$ times differentiable function $f:K\to\R$ there exist
$r$ times differentiable functions $g:\R\to\R$ and $h:\R\to\R$ such that
$f(x,y)=g(x)+h(y)$ for each point $(x,y)\in K$.

(c) We can take the graph $V_k$ of the function
$y=|x|^k$, $x\in[-1;1]$ for $k$ odd, and
$W_{k+1}=(V_{k+1}-(2,0))\cup(V_{k+1}+(2,0))$ for $k$ even.

\smallskip
{\it Proof for $k$ even.}
Let us prove that $W_{k+1}$ is $r$ times differentiably basic for each
 $0\le r\le k$. Given an $r$ times differentiable function $f:W_{k+1}\to\R$,
 take functions $h(y)=0$ and $g(x)=f(x,|x-2\sign x|^{k+1})$. 
Clearly, $h$ is $r$ times differentiable and $f(x,y)=g(x)+h(y)$ for each 
$(x,y)\in W_{k+1}$. 
Since the function $p(t)=|t|^{k+1}$ is $k$ times differentiable and $r\le k$, 
it follows that $g$ is $r$ times differentiable. 

Let us prove that $W_{k+1}$ is not $r$ times differentiably basic for
 $k$ even and each $k<r$. 
Define a function $f:W_{k+1}\to\R$ by $f(x,y)=y\sign x$.  
Clearly, $f$ is $r$ times differentiable.  
If $W_{k+1}$ is $r$ times differentiably basic, then there are
$r$ times differentiable functions $g$ and $h$ such that $f(x,y)=g(x)+h(y)$. 
For $t\in[-1;1]$ we have 
$$g(\pm2+t)+h(|t|^{k+1})=f(\pm2+t,|t|^{k+1})=\pm|t|^{k+1}.$$
Since $g$ is $(k+1)$ times differentiable and $k+1$ is odd, it follows that 
$h'(0)=+1$ and $h'(0)=-1$, which is a contradiction. 
\qed

\smallskip
{\it Proof for $k$ odd.} 
First we prove that $V_k$ is $r$ times differentiably basic for each 
$0\le r\le k$. 
Take an $r$ times differentiable function $f:V_k\to\R$.
Since $f$ is $r$ times differentiable at $(0,0)$, it follows that there exist 
$\{a_{ij}\}_{i,j=0}^r\subset\R$ such that 
$$a_{00}=f(0,0)\quad\text{and}\quad 
f(x,|x|^k)=\sum\limits_{i,j=0}^r a_{ij}x^i|x|^{kj}+o([x^2+x^{2r}]^{r/2})
\quad\text{when}\quad x\to0.$$ 
Since 
$$o([x^2+x^{2r}]^{r/2})=o_1(x^r),\quad\text{we have}\quad 
f(x,|x|^k)=a_{00}+a_{01}|x|^k+a_{10}x+\dots+a_{r0}x^r+o_2(x^r).$$ 
Take $h(y)=a_{01}y$ and $g(x)=f(x,|x|^k)-h(|x|^k)$.
Clearly, $h$ is $r$ times differentiable and $g$ is $r$ times differentiable
outside 0.  We also have $g(x)=a_{00}+a_{10}x+\dots+a_{r0}x^r+o_2(x^r)$
when $x\to0$. So $g$ is $r$ times differentiable also at 0.  

Next we prove that $V=V_1$ is not $r$ times differentiably basic for
each $1<r$. 
Define a differentiable function $f:V\to\R$ by $f(x,y)=xy$, where $y=|x|$. 
If $V$ is $r$ times differentiably basic for some $r\ge2$, then
 there are $r$ times differentiable functions 
$$g,h:\R\to\R\quad\text{such that}\quad f(x,|x|)=x|x|=g(x)+h(|x|).$$ 
Hence $g(x)-g(-x)=2x^2$ for $x\in[0;1]$. 
But this is impossible because $g$ is 2 times differentiable, hence 
for $x\to+0$ 
$$g(x)=g(0)+ax+bx^2+o(x^2)\quad\text{and}\quad g(-x)=g(0)-ax+bx^2+o(x^2).$$ 
\quad
At last we prove that $V_k$ is not $r$ times differentiably basic for
 $k$ odd and each $k<r$.  
Define a differentiable function $f:V_k\to\R$ by $f(x,y)=xy$, where $y=|x|^k$. 
If $V$ is $r$ times differentiably basic for some $r>k$, then there are
 $r$ times differentiable functions 
$$g,h:\R\to\R\quad\text{such that}\quad f(x,|x|^k)=x|x|^k=g(x)+h(|x|^k).$$ 
Hence $g(x)-g(-x)=2x^{k+1}$ for each $x\in[0;1]$. 
But this is impossible for $k$ odd because $g$ is $(k+1)$ times differentiable,
hence  for $x\to+0$ 
$$g(x)=g_0+g_1x+\dots+g_{k+1}x^{k+1}+o(x^{k+1})\quad\text{and}
\quad g(-x)=g_0-g_1x+\dots+g_{k+1}x^{k+1}+o(x^{k+1}).\quad\qed$$

\normalsize

\smallskip
\centerline{\bf References}

[Ar58] V.I. Arnold, {\it Representation of functions of some number of 
variables as superposition of functions of less number of variables (in 
Russian)}, Mat. Prosveschenie, 3 (1958), 41--61.
\linebreak
http://ilib.mirror1.mccme.ru/djvu/mp2/mp2-3.djvu?, djvuopts\&page=43

[Ar58'] V.I. Arnold, {\it Problem 6 (in Russian)}, Mat. Prosveschenie, 3 
(1958), 273-274.
\linebreak
http://ilib.mirror1.mccme.ru/djvu/mp2/mp2-3.djvu?, djvuopts\&page=243

[Ku00] V.~Kurlin, {\it Basic embeddings into products of graphs,}
Topol.Appl. 102 (2000), 113--137.

[Ku03] V. A. Kurlin, {\it Basic embeddings of graphs and the Dynnikov method of
three-pages embeddings (in Russian),} Uspekhi Mat. Nauk, 58:2 (2003), 163--164.
English transl.: Russian Math. Surveys, 58:2 (2003).

The full text of dissertation is available at 
http://maths.dur.ac.uk/$\sim$dma0vk/PhD.html

[Mi09] E. Miliczka, {\it Constructive decomposition of a function of two
variables as a sum of functions of one variable}, Proc. AMS, 137:2 (2009),
607-614.

[MK03] N. Mramor-Kosta and E. Trenklerova, {\it On basic embeddings of compacta
         into the plane,} Bull. Austral. Math. Soc. 68 (2003), 471--480.

[RZ06] D. Repov\v{s} and M. \v{Z}eljko, {\it On basic embeddings into the 
plane,} Rocky Mountain J. Math., 36:5 (2006), 1665-1677.

[Sk95] A. Skopenkov, {\it A description of continua basically embeddable in
$\R^2$,} Topol. Appl. 65 (1995), 29--48.

[St89] Y.~Sternfeld, {\it Hilbert's 13th problem and dimension,}
Lect. Notes Math. 1376 (1989), 1--49.

[Vi04] A.G. Vitushkin, \emph{Hilbert's 13th problem and related questions,}
Russian Math. Surveys,~59:1, (2004),~11--24.

[Vo81] S.M. Voronin, Funkcionalniy Analiz, 15:1 (1981), 1--17.

[Vo82] S.M. Voronin, Funkcionalniy Analiz, 16:2 (1982), 21--29.

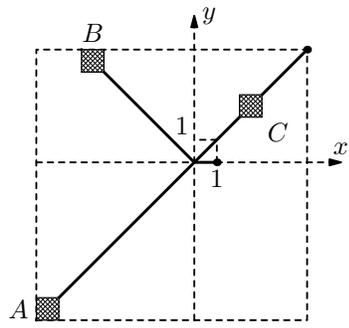
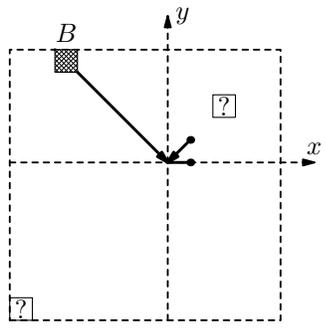
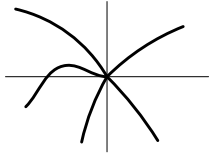
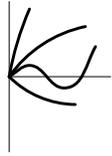
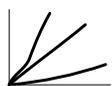
\end{document}